\documentclass{article}
\usepackage{graphicx}
\usepackage{amsmath}

\input{tcilatex}

\begin{document}

TITLE: A methodological exhibition of the theory of the identification of
Lineal Dynamic systems.

AUTHORS: Gloria Beatr\'{i}z Nu\~{n}ez Rodr\'{i}guez., Rosina Hing Cort\'{o}n
and Diosdado L. Villegas Villegas

\underline{SUMMARY}

The important role of Mathematics in the scientific development of humanity
is determined by the possibility of elaborating mathematical models of the
objects studied in the different spheres of science. In other terms, it
means, describing by means of the rigorous language of mathematics the
properties and relationships of these objects, that in turn it allows the
application of the rigorous and precise techniques of this science to the
problems related with the study of the real objects facilitating this way
the adquisition of conclusions on these and even to control their behavior.

The application of these models and methods of mathematics have been
increased significantly from the second half of the XX century, due to: the
appearance and development of the calculation techniques, the introduction
of the systemic focus in the sciences and the development of new branches
that it constitutes the Modern Applied Mathematics.

The fundamental objective of mathematics in the formation of specialists in
the different branches of the science and technical areas is to establish
that they are capable of applying the mathematical modelation and the
calculation experiment, like new technology of the scientific work in all
the spheres of the professional activity. To achieve this objective it is
necessary to develop in the students abilities for: the application of
logical forms of thinking, the application of the mathematical concepts in
the description of the objects and relationships among the real objects and
for the elaboration, selection and application of the most efficient
algorithms for the computational solution of the outlined problems.$\left[ 1%
\right] $

The identification theory and realization of the dynamic systems is a
medullary aspect in the modern control theory that consists fundamentally in
that, starting from the knowledge of the behavior entrance-exit, obtained
experimentally in the case of the identification, or given previously in the
case of the realization, to build a state model that carries out this
behavior. This content is not generally treated in the pre-graduate courses,
for the systems of multiple entrances and multiple exits.

In this work it is demonstrated that the identification theory and
realization of the lineal dynamic systems can be imparted in the technical
careers starting from the results of Mathematical Analysis and Lineal
Algebra received by the students in pre-graduate studies, without necessity
of adding new contents in the programs of this subjects

.

We propose a new form of imparting the theory of identification and
realization of the lineal dynamic systems based on the intuition and the
physical interpretation of the concepts.

\underline{INTRODUCTION}

The daily work of mathematicians necessarily consists in obtaining new
theorems, the discovery of new connections among concepts already well-known
and the introduction of new concepts. However, this is not the only
objective of mathematical investigation, the academic A. N. Kolmogorov also
considered as important objectives of mathematical investigation the
following :

1-Carry the modern logical foundations of mathematics to such a level that
they can be exposed the scholars from 14 to 15 years.

2 - To destroy the divergences among the rigorous methods of the pure
mathematicians and the non rigorous focuses of the mathematical reasonings
used by the applied mathematicians, the physicists and the technologists.$%
\left[ 2\right] $

With relationship to these objectives the mathematicians have the task of
elaborating the logical foundations of the irreproachable mathematical
intuitions, in which the rigorous calculation methods are based, arosen from
the alive intuition of the physicists and the technologists that have been
validated in the practice. But these foundations lead in occasions to
constructions so sophisticated and complex that make the mathematicians
proud, but that the physicists and engineers would not have the possibility
to assimilate if the necessary effort was not carried out to make it
accessible to them

This work meets the first objective of the investigation outlined by
Kolmolgorov and has as a result:

To develop the identification theory and realization of the lineal dynamic
systems starting from the results of Mathematical Analysis and Lineal
Algebra received by the students in pre-graduate studies with a
methodological focus that facilitates their assimilation for the students
starting from the pre-grad level.

\underline{DEVELOPMENT}

In this work a methodology is presented for the exposing of the problem of
obtaining a state model of a dynamic system, starting from the knowledge of
its entry- exit bahaviour, given by a transfer matrix, so that this is
understandable by the students from the pre-grad level whose mathematical
foundation is based on results of Functional Analysis and of the theory of
matrixes whose elements are rational functions. With this methodology it is
avoided to use sophisticated and complex mathematical theories.

The structural unit and logics of the exhibition is guaranteed by means of
the application of analysis procedures, synthesis, comparison
exemplification and of the intuition starting from the physical
interpretation of the mathematical concepts and the mathematical description
of the physical problems.

\underline{MODELLING OF SISO SYSTEMS }.

Starting from the analysis of different cases already well-known, as the
one-dimensional movement of a body under the action of an external force,
the circulation of the current in a circuit R-L, due to the application of a
voltage source in the net, studied inside the Physics, etc. will establish
that the students can arrive at the conclusion that the modelation of many
technical problems leads to the position of a system with a simple entrance
and a simple exit (SISO)

Where:

u(t): are the variables of entry of the system that contain the information
on the environment of what is disposed and the impacts in the behavior of
the study object.

y(t): are the variables of exit of the system that contain the information
which is able to be captured through the sensors.

x(t): variables of state of the system that contain the minimum quantity of
information on the system that is enough to predict their future behavior.

If the physical laws that govern the behavior of the real system are known
and also this is does not vary with regard to the adjournments in the time
and it satisfies the overlapping principle, that is to say it is stationary
and lineal, and it is also initially in rest, the mathematical pattern is
obtained:

$\left( D^{n}+a_{1}D^{n-1}+...+a_{n-1}D+a_{n}\right) \ X\ \left( t\right) \
=\ u\left( t\right) $ , \ \ \ \ X(0) = 0 \ (1)

y$\left( t\right) \ \ =\ \left( b_{0}+b_{1}D+...+b_{k}\ D^{k}\right) \
X\left( t\right) \ $, k $\leq \ n$ \ \ \ \ \ \ \ \ (2)

Through the application of Laplace Transformation , the results obtained are:

Y(S) = G(S) U(S)\ \ \ \ \ \ \ \ \ (3)

where:

G$\left( S\right) \ =\ \frac{b_{0}\ +\ b_{1\ }S\ +...+\ b_{k\ }S^{k}}{S^{n}\
+\ a_{1}\ S^{n-1}+...+\ a_{n-1}\ S\ +\ a_{n}}$\ \ \ \ \ \ \ \ \ (4)

it is denominated function of transfer of the system (which can be
interpreted as the transformation of Laplace of the answer to the unitary
impulse)

Using the transformation of coordinated , x$^{i-1}\left( t\right) =$ x$%
_{i}\left( t\right) \ $i=1,$\cdots $, n, the equations (1) and (2) can be
expressed by means of the normal system of order n

$\overset{.}{X}=AX+BU$

Y=CX\ \ \ \ \ \ \ \ \ \ \ \ \ \ \ \ \ \ \ (5)

where:

A = $\left( 
\begin{array}{cccc}
0 &  &  &  \\ 
0 &  & I_{n-1} &  \\ 
\vdots &  &  &  \\ 
-a_{n} & -a_{n-1} & \cdots & -a_{1}
\end{array}
\right) \ \ ;\ B\ =\ \left( 
\begin{array}{c}
0 \\ 
0 \\ 
\vdots \\ 
1
\end{array}
\right) \ \ y\ C\ =\ \left( 
\begin{array}{ccccccc}
b_{0} & b_{1} & \ \cdots & b_{k} & 0 & \cdots & 0
\end{array}
\right) $

being I$_{n-1}$the unitary womb of size (n -1), X $\in R^{n}$, A $\in \ $M
(n ,n),

B $\in \ $M(n,1) and C $\in $ M(1,n). The womb A it is known in the
literature like main partner of the characteristic polynomial of the
equation $(1)$

We can observe that the transfer function G(S) it is a rational function
where the degree of the polynomial of the denominator coincides with the
number of variables that intervene in the system and the polynomial of the
numerator is associated to the exit of the system .Comparing (4) and (5) it
is easy to establish a correspondence biunivoca among the pattern entry/exit
given by the transfer matrix and the state model given by the tern (A, B, C)
and pass from one to another without using the equations (1) and (2).

This type of model is denominated ''state model'' and the variable x$_{i}\ $%
; i =1,$\cdots ,$n ; are denominated state variables.

\underline{MODELLING OF MIMO} SYSTEMS

Many of the problems that are presented in practicals respond to a model of
multiple entrances and multiple exits (MIMO). We will Consider that these
systems satisfy the overlapping principle and they are invariable with
regard to the translation in time, that is to say that they are lineal and
stationary. We will also suppose that it demands the maximum degree of
pairing between the entrances and the exits. Then for the linealidad of the
system it can planted that:

y$_{i\ }=\ y_{i1}+\ y_{i2}\ +...+y_{ij\ }+...+y_{ir}\ \ \ \ \ \ \ \ \ \ \ \
\ \ \ \ \ \ \ \ \ \ \ \ \ \ \ \ \ \ \ $(6)

where y$_{ij}\ $represents the effect of the entry (j-\'{e}sima) in the exit
(i-\'{e}sima), this jeans that it is the signal of the exit (i-\'{e}sima),
if only the entry (j-\'{e}sima) was not null

~

$u_{j}\ \neq \ 0;\ u_{1}=u_{2}=...=u_{j-1}=...=\ u_{n}=0)\ \ \ \ \ \ \ \ \ $%
(7)

The same as in the case of the systems SISO where initially an analysis of
the system in rest was made, that is to say for those in u(t)$\equiv $0, y(0)%
$\equiv $0 for t$\leq $0

If we apply the Transformation of Laplace to the entrances and the exits of
(7) we can obtain the transfer function that relates the exit i-\'{e}sima
with the entrance j-\'{e}sima

G$_{ij}=\ \frac{Y_{ij}(S)}{U_{j}(S)}$ \ $\Longleftrightarrow
Y_{ij}(S)=G_{ij}\ (S)\ U_{j}(S)$ \ \ \ \ \ \ \ \ (8)

After if we apply Transformation of Laplace to (6) and we use (8) we will
obtain

Y$_{i}(S)\ =\ \overset{r}{\underset{j=1}{\sum }}G_{ij}\ (S)\ U_{j}(S);\
i=1,...,m$

whose matrical form is:

Y(S) = G(S) U(S) \ \ \ \ \ \ \ \ (9)

where:

G(S) = [G i j (S) ] ( i =1,...,m ; j =1,..., r ) is a matrice of size m x
r,Y(S) $\in $ \ R$^{m}$y U(S)$\ \in \ R^{n}\ $

Notice that equation (9) is the multidimensional analogy of equation (3)
therefore we can say (including the multidimensional and the unidimensional
case) that the Laplace transformation of the exit vector in a lineal system
is a lineal function of the Laplace transformation of the entrance vector,
defined by a matrice whose elements are functions of the variable (s) and
that it is denominated transfer matrice. The physical interpretation of this
matrice can also be made intending to apply to each entrance for shift a
very narrow and sharp pulse, of unitary area, that ideally is equal to apply
u$_{k}$(t)=$\delta _{kj}\delta $(t) in the experiment j-\'{e}simo,
k,j=1,...,r, in each experiment the exits are measured and the observations y%
$_{ij}$(t), for i=1,,m and j=1,...,r. The matrice [y$_{ij}$(t)] is the
answer matrice to the unitary impulse and the Laplace transformation [L$%
^{-1} $\{y$_{ij}$(t)\}] = G(s), is the transfer matrice.

To continue , the following problem is planted: \textquestiondown How to
obtain the transfer matrice starting from the state model?

Be a state MIMO model:

$\overset{.}{X}=AX+BU$

Y=CX \ \ \ \ \ \ \ \ \ \ \ 

Applying the Laplace Transformation, we obtain:

S X(S) = A X(S) + B U(S)

Y(S) = C X(S)

Making X (S) the subject of the formula in the first equation and
substituting in the second equation we obtain:

Y(S) = C(SI-A)$^{-1}BU(S)$

The transfer matrices is given by:

G(S) = C(SI-A)$^{-1}B$

This means that the solution of the problem outlined only requires algebraic
calculations of inversion and multiplication of matrices that can be solved
with the help of MatLab.

\underline{REALIZATION OF STATIONARY LINEAL SYSTEMS}

Let us analyze the inverse problem: Given the transfer matrice G(S) = Gi
j(S) (i=1,...,m; j=1,...,n) obtained using the formula experimentally (9)
and carrying out the measurment that were indicated when giving the
geometric interpretation of the transfer matrice.

\textquestiondown How to obtain the state pattern?

This inverse problem that one knows with the name of realization problem,
has great importance from a practical point of view and its physical
interpretation would be the following :

Given a real system, with consistent uncertainty in that we don't know the
physical laws that govern their development in time neither even define
their state variables, we only have the possibility to excite the system
with certain entrance signs and through the measuring of certain exit signs
inorder to analyze their behavior .

We will use an inductive method, that is to say, through the analysis of
particular cases we will try to infer a general methodology.

Be the system: has two entrances and an exit, this means that the
investigator can apply two types of entrance signs, for example (an electric
field and a magnetic field, or pressure and temperature) and it can measure
an exit. Therefore the transfer matrice has the type:

$[$G11(S) G12(S)]

By means of the formula (9) and the experiments described in the physical
interpretation of the transfer function, it is concluded:

1er experiment: u$_{1}$=$\delta $(t), u$_{2}$(t) =0, y(t) = Exp(-t), then G$%
_{11}\left( t\right) =\frac{Y\left( S\right) }{U_{1}\left( S\right) }=\frac{1%
}{S+1}$

2do experiment: u$_{1}$(t)= 0, u$_{2}$=$\delta $(t), y(t)=Cosh(t), then G$%
_{12}\left( t\right) =\frac{Y\left( S\right) }{U2\left( S\right) }=\frac{S}{%
S^{2}-1}$

Then

G$\left( S\right) =$ $\left[ \frac{1}{S+1}\ \frac{S}{S^{2}-1}\right] $

Comparing the results obtained in (4) and (5) we arrive for analogy to the
following submodels:

The corresponding to G$_{11}$ represented by the tern (A$_{11}$, B$_{11}$,C$%
_{11}$) where:

A$_{11}$ = (-1); B$_{11}$ = (1) and C$_{11}$ = (1)

And the corresponding to G$_{12}$ represented by the tern (A$_{12}$, B$_{12}$%
, C$_{12}$) where:

\bigskip A$_{12}=\left( 
\begin{array}{cc}
0 & 1 \\ 
1 & 0
\end{array}
\right) \ ,\ B_{12}=\left( \ 
\begin{array}{c}
0 \\ 
1
\end{array}
\right) \ \ and\ C\ \left( 
\begin{array}{ccc}
1 & 0 & 1
\end{array}
\right) $

Then the resulting system is:

$\overset{.}{x}_{1}=-x_{1}+u_{1}$

$\overset{.}{x}_{2}=x_{3}$

$\overset{.}{x}_{3}=x_{2}+u_{2}$

y$_{1}=x_{1}+x_{3}$

whose state model is defined by the tern (A, B, C) where:

A=$\left( 
\begin{array}{ccc}
-1 & 0 & 0 \\ 
0 & 0 & 1 \\ 
0 & 1 & 0
\end{array}
\right) \ \ ;\ B=\left( 
\begin{array}{cc}
1 & 0 \\ 
0 & 0 \\ 
0 & 1
\end{array}
\right) \ \ \ and\ \ \ C=\left( 
\begin{array}{ccc}
1 & 0 & 1
\end{array}
\right) $

Inorder to synthesize and to generalize the obtained methodology it will be
requested the students make an elaboration of an algorithm like the
following one:

Entry of data:

To receive G = [G ij ] (i = 1,,m ; j = 1,,n) ; \ G$_{ij}\ =\frac{%
b_{0}^{(i,j)}+\ b_{1}^{\left( i,j\right) }\ S\ +\ \cdots +b_{k_{ij}}^{\left(
i,j\right) }S^{k_{ij}}}{S^{n_{ij}}\ +\ a_{1}^{\left( i,j\right)
}S^{n_{ij}-1}+...+a_{n_{ij}}^{\left( i,j\right) }}$

Construction of sub models:

For each j = 1,...,r

For each i = 1,...,m

To construct A i j $\in $M ( n i j, n i j )

Make the last line equal to: $\ \left[ -a_{n_{ij}}^{\left( i,j\right)
}\cdots -a_{1}^{\left( i,j\right) }\right] $

Make the first column equal to: \ $\left( 
\begin{array}{c}
0 \\ 
0 \\ 
\vdots \\ 
0 \\ 
-a_{n_{ij}}^{\left( i,j\right) }
\end{array}
\right) $

Complete the matrice A i j with the unitary matrice I$_{n_{ij}-1}$

Define

\bigskip b$^{\left( i,j\right) }\in \ R^{\ n_{ij}},b^{\left( i,j\right)
}=\left( 
\begin{array}{c}
0 \\ 
0 \\ 
\vdots \\ 
1
\end{array}
\right) $

Define

C$^{\left( i,j\right) }\in \ R^{\ n_{ij}},C^{\left( i,j\right) }=\left( 
\begin{array}{ccccccc}
b_{0}^{\left( i,j\right) } & b_{1}^{\left( i,j\right) } & \cdots & 
b_{k_{ij}}^{\left( i,j\right) } & 0 & \cdots & 0
\end{array}
\right) $

Increase i

Increase j

Make A equal to the cellular diagonal matrice with cells A$_{ij}$ , (
i=1,...,m ; j = 1,...,r)

Make

B=$\left( 
\begin{array}{cccc}
b^{\left( 1,1\right) } & 0 & \cdots & 0 \\ 
b^{\left( 2,1\right) } & 0 & \cdots & 0 \\ 
\vdots & \vdots & \vdots & \vdots \\ 
b^{\left( n,1\right) } & 0 & \cdots & 0 \\ 
0 & b^{\left( 1,2\right) } & \cdots & 0 \\ 
0 & b^{\left( 2,2\right) } & \cdots & 0 \\ 
\vdots & \vdots & \vdots & \vdots \\ 
0 & 0 & \cdots & b^{\left( 1,r\right) } \\ 
0 & 0 & \cdots & b^{\left( 2,r\right) } \\ 
\vdots & \vdots & \vdots & \vdots \\ 
0 & 0 & \cdots & b^{\left( m,r\right) }
\end{array}
\right) $

make

\bigskip

C =$\left( 
\begin{array}{ccccccccccccc}
c^{\left( 1,1\right) } & 0 & \cdots & 0 & c^{\left( 1,2\right) } & 0 & \cdots
& 0 & \cdots & c^{\left( 1,r\right) } & 0 & \cdots & 0 \\ 
0 & c^{\left( 2,1\right) } & \cdots & 0 & 0 & c^{\left( 2,2\right) } & \cdots
& 0 & \cdots & 0 & c^{\left( 2,r\right) } & \cdots & 0 \\ 
\vdots & \vdots & \vdots & \vdots & \vdots & \vdots & \vdots & \vdots & 
\vdots & \vdots & \vdots & \vdots & \vdots \\ 
0 & 0 & \cdots & c^{\left( m,1\right) } & 0 & 0 & \cdots & c^{\left(
m,2\right) } & \cdots & 0 & 0 & \cdots & c^{\left( m,r\right) }
\end{array}
\right) $

This methodology will be able to be applied always and provided that they
are rational fractions, that which is presented in most of the practical
examples. If some element of G(s) it is not a rational fraction, then the
system can not be modeled with equations of the type (4).

\underline{MODEL OF MINIMAL STATE}

Until this moment we have obtained a solution of the problem for the cases
in which the elements of the transfer matrice are rational fractions
characteristic of the variable s , but can we guarantee that the solution is
unique? and in negative case ;Is that the most simple solution? We will try
to respond these questions by means of examples, avoiding the demonstration
of a mathematical theorem where it can be concluded that the solution of the
equation matricial

C ( SI- A )$^{-1}$ B = G(S) is not unique

Above all it is defined when two lineal dynamic systems are equivalent,
which is refered to as the pair

(t, x) $\in (R$, X ) where X is the state space, a phase and R x X phase
space

.

\underline{Definition:} Two lineal dynamic systems with state vectors x y x$%
\ $are algebraically equivalent whenever their numeric phase vectors are
related for all t like $\left( t,\widetilde{x}\right) $= (t, T x) where T is
a non singular matrice of size n. $\left[ 3\right] $

In the cases that we are analyzing the system of differential equations is
lineal and stationary. Be this the state model abtained

$\overset{.}{X}=AX+BU$

Y=CX \ \ \ \ \ \ \ \ (A)

Making a base change in the state space by means of the non singular
transformation T, Z=TX we obtain a system algebraically equivalent to (A)

$\overset{.}{Z}=\widetilde{A}X+\widetilde{B}U$

Y=$\widetilde{C}$ Z \ \ \ \ \ \ \ \ (B)

where

$\widetilde{A}$= T A T$^{-1}\ ,\ \widetilde{B}=TB\ ,\ \widetilde{C}=CT^{-1}$

It is demonstrated that two systems algebraically equivalent have the same
transfer matrice and it can be concluded ternas $\left( A,B,C\right) \ y\ (\ 
\widetilde{A},\widetilde{B},\widetilde{C})$ and are solutions different from
G = C ( S I - A )$^{-1}$and as a consequence this equation has infinite
solutions, since starting from any main non singular matrice T we can obtain
the corresponding ones $(\ \widetilde{A},\widetilde{B},\widetilde{C})$ .

This means that by means of the algebraic equivalence we obtain another
state model, and therefore this is not unique. For these solutions the state
space has the same dimension. KALMAN proved that there are states with
different dimensions for the same transfer matrice and it outlined the
problem of extracting the pattern of state of minimum dimension, that is to
say, given a transfer matrice G(s) to obtain the state pattern with smaller
number of variables of possible state.[3]

To obtain the minimal state it is necessary to define what we call
completely controllable system and system totally observable.

\underline{Definition}: A lineal dynamic system is completely controllable
if it is not algebraically equivalent to a system of the type:

$\overset{.}{Z}_{1}$ = A$_{11}$ z$_{1}$ + A$_{12}$ z$_{2}$+ B$_{1}$ u

$\overset{.}{Z}_{2}$= A$_{22}$ z$_{2}$

Y = c$_{1}$ z$_{1}$ + c$_{2}$ z$_{2}$

where z$_{1}$ y z$_{2}$ are vectors of n$_{1}$ y n$_{2}$ = n - n$_{1}$
components respectively

In other words, if the lineal system (A) is algebraically equivalent to a
system that has for state model terna $(\ \widetilde{A},\widetilde{B},%
\widetilde{C})$ such that:

$\widetilde{A}$=$\left( 
\begin{array}{cc}
A_{11} & A_{12} \\ 
0 & A_{22}
\end{array}
\right) \ ;\ \widetilde{B}=\left( 
\begin{array}{c}
B_{1} \\ 
0
\end{array}
\right) \ \ y\ \widetilde{C}$\ \ \ whatever

The system (A) is not completely controllable, that is to say, it is
possible to find a system of coordinates in which the variable z$_{i}$ is
separated in two groups z$_{1}=\left( z^{1},z^{2},\cdots ,z^{n_{1}}\right) \
,\ \ y\ \ z_{2}\left( z^{n_{1}+1},\cdots ,z^{n}\right) $, which are
denominated controllable and not controllable respectively.

The application of the definition is annoying since if the initial system
has the form of $(\ \widetilde{A},\widetilde{B},\widetilde{C})$ it is
evident that it is controllable, but if it doesn't have, we can not affirm
nothing since there are infinite possibilities of algebraic transformations.
For such a reason we will apply without demonstration a characterization
based on the concept of range of a matrice. We could have used the variant
of defining the concept however by means of the characterization that has
possible algorithms, it is a very abstract definition. The form that we use
gives the possibility to give an intuitive idea, a physical interpretation
of the concept and for this reason we have selected it.

\bigskip

The caracterization is the following:

In the system:

$\overset{.}{X}=AX+BU$

Y=CX \ \ \ \ \ \ \ \ \ \ \ \ \ \ \ \ \ \ (A) \ where A$\in $M (n, n )

We consider that B is a matrice of the size n x r and M c = [ B AB$\cdots $ A%
$^{n-1}$ B ] a matrice of n rows and n x r coloums called controllable
matrice.

\bigskip

''The condition necesary and sufficient for the system (A) to be totally
controllableis that

R(M $_{c}$) = n '' [4].

To obtain the R(M c) an algorithm is deduced to find the k ( k $\leq $\ n )
columns lineally independent of M$_{c}$ without working jointly with M $_{c}$
columns.

To reciever A y B

To define H, a group of at most n vectors R$^{n}$, inicially empty

for $\alpha $= 1 to $\alpha $= r choose b$_{\alpha }\in \ $B

make j = 0

while r [ b $_{\alpha }$\ \ Ab$_{\alpha }$ \ A$^{j}b_{\alpha }$] = j+1

Include the vector coloumn A$^{j}b_{\alpha }$ in H

Increase j

Increase $\alpha $

Return H

End

\underline{Definition}: The space generated by all the columns lineally
independent of the matrice M $_{c}$ is called controllable space and we will
denote it as E$_{c}$.

\underline{Property}: E $_{c}$ is an invariant [5] with relationship to the
system

$\overset{.}{X}=AX+BU$

Y=CX

\ \ \ \ \ \ \ \ \ \ \ \ 

\textquestiondown How to separate the completely controllable and not
completely controllable spaces?

If E $_{c}$ = \TEXTsymbol{<} f$_{1}\cdots $ f $_{k}$ \TEXTsymbol{>} being \{
f$_{1}\cdots $ f $_{k}\}$ a base of E$_{c}$ and let us complete with the
canonical vectors \{ e$_{1}$,$\cdots $ e$_{n}$ \} a base of En. If the
lineal transformation T=[f$_{1}\cdots f_{k}\ e_{k+1}\cdots e_{n}]\ $of E n
such that

Z = T X = x$_{1}$ f$_{1}$ +$\cdots $ + x $_{k}$ f $_{k}$+ x $_{k+1}$ e$%
_{k+1} $ +$\cdots $ + x$_{n}$ e $_{n}$

.

It is demonstrated that by means of the matrice T a system is obtained
algebraically equivalent to (A) whose state model is the terna $(\ 
\widetilde{A},\widetilde{B},\widetilde{C})$ where:

$\widetilde{A}$=$\left( 
\begin{array}{cc}
A_{11} & A_{12} \\ 
0 & A_{22}
\end{array}
\right) \ ;\ \widetilde{B}=\left( 
\begin{array}{c}
B_{1} \\ 
0
\end{array}
\right) \ \ y\ \widetilde{C}\ $whatever.

It means,

$\overset{.}{Z}_{1}$ = A$_{11}$ z$_{1}$ + A$_{12}$ z$_{2}$+ B$_{1}$ u

$\overset{.}{Z}_{2}$= A$_{22}$ z$_{2}$

Y = $\widetilde{C}z$

We will now analyze the dual concept with regard to the controllable concept.

\underline{Definition}: A lineal dynamic system is totally observable if
this it is not algebraically equivalent to a system of the type:

$\overset{.}{Z}_{1}$ = A$_{11}$ z$_{1}$ + B$_{1}$ u

$\overset{.}{Z}_{2}$= A$_{21}$ z$_{1}+$A$_{22}$ z$_{2}$ + B$_{2}$ u

Y =c$_{1}z_{1}$

where z$_{1}$ y z$_{2}$ are vectors of n$_{1}$and n$_{2}$ = n- n$_{1}$
components respectively. [3]

In other words, if the lineal system (A) is algebraically equivalent to a
system that has for state model the tern where:

\bigskip $\widetilde{A}$=$\left( 
\begin{array}{cc}
A_{11} & 0 \\ 
A_{21} & A_{22}
\end{array}
\right) \ ;\widetilde{C}=\left( c_{1}\ \ \ 0\right) \ $ and $\widetilde{B}$
whatever

The system (A) it is not totally observable, that is to say, it is possible
to find a system of coordinates in which the variable z$_{i}$ is separated
in two groups z$_{1}=\left( z^{1},z^{2},\cdots ,z^{n_{1}}\right) \ ,and\ \
z_{2}\left( z^{n_{1}+1},\cdots ,z^{n}\right) $ which are denominated
observables and non observables respectively.

The same as in the study of the control the application of the definition is
annoying and for such a reason we apply the following characterization:

In the system:

$\overset{.}{X}=AX+BU$

Y=CX \ \ \ \ \ \ \ \ \ \ \ \ \ \ \ \ \ \ (A) where A $\in $ M (n,n )

Consider that C is a matrice of the size m x n and

M$_{o}=\left( 
\begin{array}{c}
C \\ 
CA \\ 
\vdots \\ 
CA^{n-1}
\end{array}
\right) $

A matrice of n m rows and n coloumns called matrice of observability

''The condition necesary and suficient for the system (A) to be totally
observable is R (M $_{o}$) = n '' [4]

As transposing a matrice its range don't vary, it have:

\bigskip r(M$_{o})=r\left( 
\begin{array}{cccc}
C^{t} & A^{t}C^{t} & \cdots & (A^{t})^{n-1}C^{t}
\end{array}
\right) $

Note: The system:

$\overset{.}{X}=AX+BU$

Y=CX \ \ \ \ \ \ 

It is totally observable if and alony if X= A$^{t}$ X+ C$^{t}$U is
completely controllable, in this sense it is to say that both concepts are
dual. Consequently the R (M$_{o}$) is calculated equally to R( M $_{C}$ ).

In this case the space generated by the lineally independent rows of Mo is
denominated observable space and we will denote it E o and their complement
ortogonal is denominated non observable space (E no )

How to seperate the sub-spaces E no y E o ?

Let us suppose that E$_{C}\ =\left[ g_{1}\cdots g_{n-k_{1}}\right] \ $and E$%
_{NO}=\left[ v_{1}\cdots v_{k_{1}}\right] $ , if we complete the base in the
non observabile sub-space to a base in the whole space we obtain a matrice T$%
_{1}$.

If the lineal transformation T$_{1}$=[v$_{1}\cdots v_{k}\ e_{k+1}\cdots
e_{n}]\ $of E n such that

Z = T$_{1}$ X = v$_{1}$x$_{1}$ +$\cdots $ + v$_{k_{1}}$x $_{k_{1}}$+ v $%
_{k_{1}+1}$ e$_{k_{1}+1}$ +$\cdots $ + v$_{n}$ e $_{n}$.

Let us demonstrate that by means of the matrice T$_{1}$ a system is obtained
that is algebraically equivalent to (A) whose state model is the tern $(\ 
\widetilde{A},\widetilde{B},\widetilde{C})$. where:

\bigskip $\widetilde{A}$=$\left( 
\begin{array}{cc}
A_{11} & A_{12} \\ 
0 & A_{22}
\end{array}
\right) \ $; $\widetilde{B}$ whatever and $\widetilde{C}=[0\ C_{2}]\ $, or
if,

$\overset{.}{Z}_{1}$ = A$_{11}$ z$_{1}$ + A$_{12}$ z$_{2}$+ B$_{1}$ u

$\overset{.}{Z}_{2}$ = A$_{22}$ z$_{2}$+ B$_{2}$ u

Y = C$_{2}$ z$_{2}$

Demostration:

In effect: E$_{no}$ = \{ z : z$_{2}$ = 0 \} where Z =$\left( 
\begin{array}{c}
z_{1} \\ 
z_{2}
\end{array}
\right) $

If $\widetilde{A}$=$\left( 
\begin{array}{cc}
A_{11} & A_{12} \\ 
A_{21} & A_{22}
\end{array}
\right) \ ;\ \widetilde{B}=\left( 
\begin{array}{c}
B_{1} \\ 
B_{2}
\end{array}
\right) \ and\ C=[c_{1}\ c_{2}]$ , therefore the system is:

$\overset{.}{Z}_{1}$ = A$_{11}$ z$_{1}$ + A$_{12}$ z$_{2}$

$\overset{.}{Z}_{2}$ = A $_{21}$ z $_{1}$+A$_{22}$ z$_{2}$

Y = c$_{1}$ z$_{1}$ + c$_{2}$ z$_{2}$

We choose z$_{0}$ $\in \ $E$_{no}$ where Z$_{0}$ =$\left( 
\begin{array}{c}
z_{10} \\ 
0
\end{array}
\right) $therefore z(t) $\in $E$_{C}$; Z$(t)$ =$\left( 
\begin{array}{c}
z_{1} \\ 
0
\end{array}
\right) $ belong to E$_{no}$ then, c$_{1}$ z$_{10}$ =0

Making z equal to the vectors canonicals of E $_{no}$ in z$_{10}^{1}c_{1}+$z$%
_{10}^{2}c_{2}+\cdots +$z$_{10}^{k}c_{k_{1}}=0$

we obtain:

c$_{1}=c_{2}=\cdots =c_{k_{1}}=0$

and $\widetilde{C}=[0\ C_{2}]$

As z$_{2}=$ 0$\ \Rightarrow $\ z$_{2}$= 0$\ \Rightarrow $\ A$_{21}$ z$_{10}$
= 0 $\Rightarrow $\ A$_{21}$ = 0, as z$_{10}\ $can be whatever arbituary in
the sub-space of non-controllability. Then

$\widetilde{A}$=$\left( 
\begin{array}{cc}
A_{11} & A_{12} \\ 
0 & A_{22}
\end{array}
\right) \ $obtained the desired

Using analysis procedures and synthesis we can outline the following
methodology:

\underline{ALGORITHM FOR THE CANONICAL DECOMPOSITION OF A SISTEMA MIMO,
LINEAL AND STATIONARY}

Inicial data: A $\in \ $M (n, n), B $\in $ M (n, r), C $\in $ M (m, n).

Start

1- Calculate: rank [ B AB $\cdots $A$^{n-1}$ B ] = k and extract a base from
the sub-space E c generated by the coloumns of this matrice:

E c = $\langle e_{1}\cdots e_{k}\rangle $ , dim E c = k $\preceq \ $n, E $%
_{C}$ $\subset $ n $^{n}$,E $\in $ M (n, k),

E = [ e$_{1}\cdots $ e $_{k}$ ]

2- Calculate rank$\left( 
\begin{array}{cccc}
C^{t} & A^{t}C^{t} & \cdots & (A^{t})^{n-1}C^{t}
\end{array}
\right) $ = p and extract a base from the sub-space E o generated by the
columns of this matrice:

E o =$\langle f_{1}\cdots f_{p}\rangle $ , dim E o = p $\leq $\ n, E$%
_{o}\subset \ $E $^{n}$, F $\in \ $M (n, p), F = [ f$_{1}\cdots $f p ]

3- If E $_{no}$ = E $_{o}^{\perp }$. calculate a base from the sub-space E c 
$\cap \ $E $_{no}$, finding a fundamental system of solutions from the
system:

F$^{t}$ E $\alpha $ = 0, $\alpha \in \ $R$^{k}$,

being: E$_{c}$ = \{ v$\in $ E $^{n}$ : v = E $\alpha ,\forall \alpha \in \
R^{k}$ \}

and E $_{no}$ = \{ v $\in \ $E$^{n}$ : F$^{t}$v = 0 \}

Being the fundamental system of solutions lineally idependent: \{$\alpha
_{1} $ $\alpha _{2}\cdots \alpha _{q}$\}, where: q $\leq \ $min ( k, n p ) ,

dim E $_{no}$ = n - p = l

g i = $\sum \alpha _{ij}\ e_{j},\ $i=1,$\cdots ,q,$ E c$\ \cap $ E$_{no}$ =$%
\langle \ g_{1},\cdots ,g_{k}\rangle $ = E cno

4- To complete the vectorial g i , i = 1,$\cdots ,$q until a base in E c
applying the method of elimination of Gauss of the matrice [ g$_{1}$ g$%
_{2}\cdots $g$_{q}$ e$_{1}$e$_{2}\cdots $ e$_{k}$ ]

Note: When applying the elimination method all the vectors

g i , i = 1,$\cdots $, q. were used to build the escalar matrice. The
vectors e$_{i_{j}}$

j = 1,$\cdots $ ,k -q that are used in this submatriz is those that complete
the base.

5- To find a base in E $no$ =$\langle h_{1}\ h_{2}\ \cdots h_{l}\rangle $

6- To complete the vectorial g i , i = 1, $\cdots ,$q until a base in E no,
applying the one

method of elimination of Gauss of the matrice [ g$_{1}$ g$_{2}\cdots $g$_{q} 
$ h$_{1}$h$_{2\cdots }$h$_{l}$]

Note: All the vectorial g i , i = 1,$\cdots $,q was used, that form the main
step, the vectors,h$_{i_{j}}$ , j = 1,$\cdots $, k-q that are used to
complete this matrice is those that complete the base.

7- To complete until a base of E$^{n}$if it is necessary ( k + l - q 
\TEXTsymbol{<} n ) applying the elimination method of Gauss to the matrice [
g$_{1}$ g$_{2}\cdots $g$_{q}$ e$_{i_{1}}e_{i_{2}}\cdots $ e$_{i_{k-q}}\
h_{j_{1}}\ h_{j_{2}}\cdots h_{j_{l-q}}\ v_{1}\ v_{2}\cdots v_{n}$ ]

where the v i , i = 1,$\cdots $,n are canonical vectors. To form the
submatriz step only the v$_{i}$ were used. When it is not possible to use
the remaining columns these v$_{i}$ used are those that complete the base.

8- Construct the matrice of change of base:

If k + l - q = n then T =[ g$_{1}$ g$_{2}\cdots $g$_{q}$ e$%
_{i_{1}}e_{i_{2}}\cdots $ e$_{i_{k-q}}\ h_{j_{1}}\ h_{j_{2}}\cdots
h_{j_{l-q}}]$

If k + l - q \TEXTsymbol{<} n then \ T=[ g$_{1}$ g$_{2}\cdots $g$_{q}$ e$%
_{i_{1}}e_{i_{2}}\cdots $ e$_{i_{k-q}}\ h_{j_{1}}\ h_{j_{2}}\cdots
h_{j_{l-q}}\ v_{m_{1}}\ v_{m_{2}}\cdots v_{m_{i-q}}]$

9-To make the reduction to the canonical form finding

$\widetilde{A}$=T$^{-1}$AT, $\widetilde{B}=T\ ^{-1}B\ \ y\ \widetilde{C}=CT$

The state variables with the index \{ i$_{1}\ $i$_{2}\cdots $ i$_{k-q}$ \}
are completely controllable and totally observables and they constitute the
states of minimal representation of the analyzed system.

In the canonical representation, we divide the state variables in 4 groups:

Controllable and non observables, controllable and observables, not
controllable and non observables, not controllable and observables.

To determine the efficiency of this algorithm you can analyze their
temporary complexity. The number of operations is of order K$^{3}$, where K
= max (n r,n m), since in all the procedures you can apply the method of
elimination of Gauss and the number of operations in this is of the order of
the third power of the maximum between the number of rows and the number of
columns [6], for what its temporary complexity is of the order of K$^{\frac{3%
}{2}}$, it is therefore an algorithm polynomial.

Notice you the only jump in the logical structure of this exhibition takes
place when they are accepted without demonstration, the practical approaches
of control and observation total, really it is enough to demonstrate only
one of these criterion, for the relationship of duality that exists among
them. These criterions are of great practical utility for their algorithm
form however they have an abstract character and the demonstrations are
sophisticated points.

The demonstration of the necessity of the condition on the range of the
control matrice can be demonstrated in a relatively simple way, if we lean
on in the Theorem of Hamilton-Cayle [7] that in spite of being a very
important theorem of the algebra, and specially useful inside the Theory of
Control, and it is not included in the undergraduate courses. Several
variants of demonstration of this Theorem appear in [8], among them a very
simple and elegant one based on purely algebraic arguments.

The demonstration of the sufficiency of the criterion is more complex, a
relatively simple variant can be found in [7] based on the concept of group
of accessibility and using the theory of the differential equations.

In [9] the necessity is demonstrated using Hamilton Cayley's Theorem and the
formula of Sylvester that it is a generalization of the formula of
interpolation of Lagrange for the case of the matrices. These information
should be offered to the students, so that the interested ones in deepening
in this theory have the appropriate bibliography.

\underline{CONCLUSIONS}

It is developed an exhibition of the position and solution of the
realization problem or identification of the lineal and stationary dynamic
systems, without making use of the results of complex mathematical theories,
related with the theory of the moments, the representation theorems of
functional and lineal operators in spaces topologies and the theory of the
matrices with elements in the field of the forms models - univaried that
serve as foundation mathematics for the solution of this problem.

The unit of the exhibition of the realization theory or identification is
achieved through the intuition, based on the interpretation of such physical
concepts as: overlapping, control and observation, the exemplification, the
analysis, the synthesis and the comparison.

A solution of the realization problem for the systems SISO is obtained in
immediate form starting from the analysis and the comparison of the position
and of the solution of the direct problem: '' To find the transfer function
applying the definition from this to the state pattern.''

For study of the systems MIMO study, one proceeds to decompose these, by
virtue of the lineality, in a subsystem group SISO. each one of which are
applied the previous results, obtaining a tena (A,B,C) in which thematrice A
is a cellular diagonal matrice. The generalization of this result is
expressed by means of the elaboration of an algorithm that allows to obtain
(A,B,C) for any particular case.

The existence of an infinite number of solutions of the equation matricial

C (S I - A )$^{-1}$B = G(S)

puts on evidence in the fact that any system algebraically equivalent to the
one obtained, also carries out the behavior entrance-exit defined by the
transfer womb G(s).

Starting from the theory of the subspaces of a vectorial space, the systems
algebraically equivalent and the physical concepts of total observation and
total control are based in the obtaining of the canonical decomposition of
the obtained pattern and is a completely controllable subsystem and totally
observable and carries out the behavior pattern entrance-exit of G(s).

The fact that the subsystem of the state variables totally observables and
completely controllable it determines the realization minimal of G(s), it is
easily accepted by the students, for the evident reduction of the dimension
of the state space and the physical interpretation of the control concepts
and total observation.

\underline{BIBLIOGRAPHY}

$\left[ 1\right] $ Hing Cort\'{o}n Rosina: Programs for the Applied
development Mathematicas in the UCLV, work presented in the event
international Pedagogy 97.

$\left[ 2\right] $ Kolmogorov A. N,: The Mathematical in the Modern World,
Izvestia,1962.

$\left[ 3\right] $ Kalman R.E.: Mathematical Description of Linear Dynamical
Systems, J.S.I.A.M. Control, Ser.A, Vol. 1, N o 2, U.S.A., 1963.

$\left[ 4\right] $Sontang E.: Mathematical Control Theory, Deterministic
Finite DimensionalSystems, Springer-Verlag Irc., New York, 1990.

$\left[ 5\right] $ Tijonov A.N, Kostomarov D.P.: Something about Applied
Mathematics, Editorial MIR, Moscow, 1983.

$\left[ 6\right] $ Fadeev Faddeva: Computational Methods of Linear Algebra,
Revolutionary Editions, Cuba, 1963

$\left[ 7\right] $ Alexandrov V.V.: Mathematical Modelaci\'{o}n of Dynamic
Systems, U.H., Cuba 1988.

$\left[ 8\right] $ Bellman R.: Introduction to the An\'{a}lisis Matricial,
Editorial Revert\'{e}, Barcelona, 1965

.

$\left[ 9\right] $ Ogata Katsuhiko: Engineering of Modern Control, Prentice
Hull Hispanoamerican S. A., Mexico, 1993.

\end{document}